\documentclass[11pt]{article}
\usepackage{amsmath, amsthm, amscd, amsfonts, amssymb, graphicx, color,cite,mathpazo}
\usepackage[sans]{dsfont}
\usepackage[bookmarksnumbered, colorlinks, plainpages]{hyperref}
\usepackage{bbm}
\hypersetup{citecolor={blue}}

\numberwithin{equation}{section}
\newtheorem{theorem}{Theorem}[section]
\newtheorem{lemma}[theorem]{Lemma}

\theoremstyle{definition}
\newtheorem{definition}[theorem]{Definition}

\newtheorem{remark}[theorem]{Remark}
\newtheorem{proposition}[theorem]{Proposition}
\newtheorem{corollary}[theorem]{Corollaire}

\begin{document}
	\thispagestyle{empty} \setcounter{page}{1}

	\begin{center}
		{\LARGE\bf Midpoint and trapezoid type inequalities for multiplicatively convex
functions}

		\vskip.40in
		
		{\bf  A. Berhail$^{1}$, B. Meftah$^{2}$ } \\[2mm]
		
		{\footnotesize
			
			$^{1}, ^{2}$ Department of Mathematics, 08 May 1945 University-Guelma, Algeria }

	\end{center}
	\vskip 5mm
	
	\noindent {\footnotesize{\bf Abstract.}~ {In this paper, we first prove two new identities for multiplicative
differentiable functions. Based on this identity, we establish a midpoint
and trapezoid type inequalities for multiplicatively convex functions.
Applications to special means are also given.}}

	\vskip .15in
	{
		\footnotetext{\textbf{$^\S$ Corresponding author}: 
		Berhail Amel (E-mail: berhail.amel@univ-guelma.dz)}
	    \footnotetext{Amel Berhail( E-mail: berhail.amel@univ-guelma.dz)}
		\footnotetext{Badri Meftah (E-mail:badrimeftah@yahoo.fr)
	}
}

	\baselineskip=14pt

	\setcounter{section}{0} \numberwithin{equation}{section}

\section{\protect\bigskip \textbf{Introduction}}

The concept of convexity plays an important and very central role in many
areas, and has a close relationship in the development of the theory of
inequalities, which is an important tool in the study of some properties of
solutions of differential quations as well as in the error estimates of
quadrature formulas.
Let $I$ be an interval of real numbers

\begin{definition}[{$[ 17] $}]
A function $f:I\rightarrow \mathbb{R}$ \ is said to be convex, if \ \ 
\begin{equation*}
f\left( tx+\left( 1-t\right) y\right) \leq tf\left( x\right) +\left(
1-t\right) f(y)
\end{equation*}%
holds for all $x,y\in I$ and all $t\in \lbrack 0,1]$.
\end{definition}

\begin{definition}[{$\left[ 17\right] $}]
A positive function $f:I\rightarrow \mathbb{R}$ is said to be
logarithmically convex or multiplicatively convex , if 
\begin{equation*}
f(tx+(1-t)y)\ \leq \left[ f(x)\right] ^{t}\left[ f(y)\right] ^{1-t}
\end{equation*}%
holds for all $x,y\in I$ and all $t\in \lbrack 0,1]$.
\end{definition}

The fundamental inequality for convex functions is undoubtedly the
Hermite-Hadamard inequality, which can be stated as follows: For every
convex function $f$ on the interval $\left[ a,b\right] $ with $a<b$, we have%
\begin{equation}
f( \tfrac{a+b}{2}) \leq \tfrac{1}{b-a}\int\limits_{a}^{b}f(
x) dx\leq \tfrac{f( a) +f( b) }{2}.
\end{equation}%
If the function $f$ is concave, then (1.1) holds in the reverse direction
see $\left[ 17\right] $.

In $\left[ 6\right] $, Dragomir and Agarwal established some trapezoid type
inequalities for functions whose absolute value of the first derivatives are
convex%
\begin{equation}
\left\vert \tfrac{f\left( a\right) +f\left( b\right) }{2}-\tfrac{1}{b-a}%
\int\limits_{a}^{b}f\left( x\right) dx\right\vert \leq \tfrac{b-a}{8}\left(
\left\vert f^{\prime }\left( a\right) \right\vert +\left\vert f^{\prime
}\left( b\right) \right\vert \right) .
\end{equation}

In $\left[ 16\right] $, Pearce and Pe\v{c}ari\'{c} proved some midpoint type
inequalities for functions whose absolute value of the first derivatives are
convex 
\begin{equation}
\left\vert f\left( \tfrac{a+b}{2}\right) -\tfrac{1}{b-a}\int%
\limits_{a}^{b}f\left( x\right) dx\right\vert \leq \tfrac{b-a}{8}\left(
\left\vert f^{\prime }\left( a\right) \right\vert +\left\vert f^{\prime
}\left( b\right) \right\vert \right) .
\end{equation}%
Concerning some papers dealing with some quadrature see $\left[ 9,10,11%
\right] $ and references therein.

In $1967$, Grossman and Katz, created the first non-Newtonian calculation
system, called geometric calculation. Over the next few years they had
created an infinite family of non-Newtonian calculi, thus modifying the
classical calculus introduced by Newton and Leibniz in the $17$th century
each of which differed markedly from the classical calculus of Newton and
Leibniz known today as the non-Newtonian calculus or the multiplicative
calculus, where the ordinary product and ratio are used respectively as sum
and exponential difference over the domain of positive real numbers see $%
\left[ 8\right] $. This calculation is useful for dealing with exponentially
varying functions.

The complete mathematical description of multiplicative calculus was given
by Bashirov et al. $\left[ 4\right] $. Also in the literature, there remains
a trace of a similar calculation proposed by the mathematical biologists
Volterra and Hostinsky $\left[ 18\right] $ in $1938$ called the Volterra
calculation which is identified as a particular case of multiplicative
calculation.

Recently, Ali et al. $\left[ 1\right] $, gave the analogue of the
Hermite-Hadamard inequality for multiplicatively convex functions

\begin{theorem}
Let $f$ be a positive and multiplicatively convex function on interval $%
[a,b] $, then the following inequalities hold%
\begin{equation}
f\left( \tfrac{a+b}{2}\right) \leq \left( \int\limits_{a}^{b}f\left(
x\right) ^{dx}\right) ^{\frac{1}{b-a}}\leq \sqrt{f\left( a\right) f\left(
b\right) }.
\end{equation}
\end{theorem}

Also, they proved the inequalities for the product and the quotient of two
multiplicatively convex functions. In $\left[ 2\right] $, Ali et al. studied
the Hermite-Hadamard type inequalities for multiplicatively $\phi $-convex
and log-$\phi $-convex functions. \\ In $\left[ 12\right] $ and $\left[ 14%
\right] $, \"{O}zcan generalized the results obtained in $\left[ 1\right] $
for $s$-convex and $h$-convex functions respectively.\\ In $\left[ 13\right] $%
, \"{O}zcan gave the analogue of Hermite-Hadamard type inequalities for
multiplicatively peinvex functions. In $\left[ 15\right] $, \"{O}zcan
proposed the Hermite-Hadamard type inequalities for multiplicatively $h$%
-peinvex functions. \\
In $\left[ 3\right] $ Ali et al. investigate the
Ostrowski and the Simpson type inequalities for multiplicatively convex
functions.

Motivated by all the above papers, in this study we propose two new
identities for multiplicatively differentiable functions, based on these
identities we establish a midpoint and trapezoid type inequalities for
multiplicatively convex functions. Applications to special means are also
given.

\section{\textbf{Preliminaries}}

In this section we begin by recalling some definitions, properties and
notions of derivation as well as multiplicative integration

\begin{definition}[{$\left[ 4\right] $}]
Let $f:%
\mathbb{R}
\rightarrow 
\mathbb{R}
^{+}$ be a positive function. The multiplicative derivative of the function $%
f$ noted by $f^{\ast }$ is defined as follows%
\begin{equation*}
\tfrac{d^{\ast }f}{dt}=f^{\ast }\left( t\right) =\underset{h\rightarrow 0}{%
\lim }\left( \tfrac{f\left( t+h\right) }{f\left( t\right) }\right) ^{\frac{1%
}{h}}.
\end{equation*}
\end{definition}

\begin{remark}
If $f$ has positive values and is differentiable at $t$, then $f^{\ast }$
exists and the relation between $f^{\ast }$ and ordinary derivative $%
f^{\prime }$ is as follows:%
\begin{equation*}
f^{\ast }\left( t\right) =e^{\left( \ln f\left( t\right) \right) ^{\prime
}}=e^{\frac{f^{\prime }\left( t\right) }{f\left( t\right) }}.
\end{equation*}
\end{remark}

The multiplicative derivative admits the following properties:

\begin{theorem}[{$\left[ 4\right] $}]
Let $f$ and $g$ be multiplicatively differentiable functions, and $c$ is
arbitrary constant. Then functions $cf,fg,f+g,f/g$ and$fg$ are $^{\ast }$
differentiable and
\end{theorem}

\begin{itemize}
\item $\left( cf\right) ^{\ast }\left( t\right) =f^{\ast }\left( t\right) .$

\item $\left( fg\right) ^{\ast }\left( t\right) =f^{\ast }\left( t\right)
g^{\ast }\left( t\right) .$

\item $\left( f+g\right) ^{\ast }\left( t\right) =f^{\ast }\left( t\right) ^{%
\frac{f\left( t\right) }{f\left( t\right) +g\left( t\right) }}g^{\ast
}\left( t\right) ^{\frac{g\left( t\right) }{f\left( t\right) +g\left(
t\right) }}.$

\item $\left( \frac{f}{g}\right) ^{\ast }\left( t\right) =\frac{f^{\ast
}\left( t\right) }{g^{\ast }\left( t\right) }.$

\item $\left( f^{g}\right) ^{\ast }\left( t\right) =f^{\ast }\left( t\right)
^{g\left( t\right) }f\left( t\right) ^{g^{\prime }\left( t\right) }.$
\end{itemize}

In $\left[ 4\right] $, Bashirov et al. introduced the concept of the $^{\ast
}$ integral called multiplicative integral which is noted $\underset{a}{%
\overset{b}{\int }}\left( f\left( t\right) \right) ^{dt}$

Which the sum of the terms of the product is used in the definition of a
classical Riemann integral of $f$ over $[a,b]$, while the product of the
terms raised to the power is used in the definition of the multiplicative
integral of $f$ over $[a,b]$.

The relationship between the Riemann integral and the multiplicative
integral is as follows:

\begin{proposition}[ 4]
If $f$ is Riemann integrable on $[a,b]$, then $f$ is multiplicative
integrable on $[a,b]$ and%
\begin{equation*}
\underset{a}{\overset{b}{\int }}\left( f\left( t\right) \right) ^{dt}=\exp
\left( \underset{a}{\overset{b}{\int }}\ln \left( f\left( t\right) \right)
dt\right) .
\end{equation*}
\end{proposition}

Moreover, Bashirov et al. show that multiplicative integral has the
following results and properties:

\begin{theorem}[{$\left[ 4\right] $}]
Let $f$ be a positive and Riemann integrable on $[a,b]$, then $f$ is
multiplicative integrable on $[a,b]$ and
\end{theorem}

\begin{itemize}
\item $\underset{a}{\overset{b}{\int }}\left( \left( f\left( t\right)
\right) ^{p}\right) ^{dt}=\left( \underset{a}{\overset{b}{\int }}\left(
f\left( t\right) \right) ^{dt}\right) ^{p}.$

\item $\underset{a}{\overset{b}{\int }}\left( f\left( t\right) g\left(
t\right) \right) ^{dt}=\underset{a}{\overset{b}{\int }}\left( f\left(
t\right) \right) ^{dt}\underset{a}{\overset{b}{\int }}\left( g\left(
t\right) \right) ^{dt}.$

\item $\underset{a}{\overset{b}{\int }}\left( \frac{f\left( t\right) }{%
g\left( t\right) }\right) ^{dt}=\frac{\underset{a}{\overset{b}{\int }}\left(
f\left( t\right) \right) ^{dt}}{\underset{a}{\overset{b}{\int }}\left(
g\left( t\right) \right) ^{dt}}.$

\item $\underset{a}{\overset{b}{\int }}\left( f\left( t\right) \right) ^{dt}=%
\underset{a}{\overset{c}{\int }}\left( f\left( t\right) \right) ^{dt}%
\underset{c}{\overset{b}{\int }}\left( f\left( t\right) \right) ^{dt},a<c<b.$

\item $\underset{a}{\overset{a}{\int }}\left( f\left( t\right) \right)
^{dt}=1$ and $\underset{a}{\overset{b}{\int }}\left( f\left( t\right)
\right) ^{dt}=\left( \underset{b}{\overset{a}{\int }}\left( f\left( t\right)
\right) ^{dt}\right) ^{-1}.$
\end{itemize}

\begin{theorem}[{Multiplicative Integration by Parts $\left[ 4\right] $}]
Let $f:[a,b]\rightarrow 
\mathbb{R}
$ be multiplicative differentiable, let $g:[a,b]\rightarrow 
\mathbb{R}
$ be differentiable so the function $f^{g}$ is multiplicative integrable, and%
\begin{equation*}
\underset{a}{\overset{b}{\int }}\left( f^{\ast }\left( t\right) ^{g\left(
t\right) }\right) ^{dt}=\tfrac{f\left( b\right) ^{g\left( b\right) }}{%
f\left( a\right) ^{g\left( a\right) }}\times \tfrac{1}{\underset{a}{\overset{%
b}{\int }}\left( f\left( t\right) ^{g^{\prime }\left( t\right) }\right) ^{dt}%
}.
\end{equation*}
\end{theorem}

\begin{lemma}[{$\left[ 3\right] $}]
Let $f:[a,b]\rightarrow 
\mathbb{R}
$ be multiplicative differentiable, let $g:[a,b]\rightarrow 
\mathbb{R}
$ and let $g:J\subset 
\mathbb{R}
\rightarrow 
\mathbb{R}
$ be two differentiable functions. Then we have%
\begin{equation*}
\underset{a}{\overset{b}{\int }}\left( f^{\ast }\left( h\left( t\right)
\right) ^{h^{\prime }\left( t\right) g\left( t\right) }\right) ^{dt}=\tfrac{%
f\left( b\right) ^{g\left( b\right) }}{f\left( a\right) ^{g\left( a\right) }}%
\times \tfrac{1}{\underset{a}{\overset{b}{\int }}\left( f\left( h\left(
t\right) \right) ^{g^{\prime }\left( t\right) }\right) ^{dt}}.
\end{equation*}
\end{lemma}

\section{\textbf{Main results}}

In order to prove our results, we need the following lemmas

\begin{lemma}\label{lem2}
Let $f:$ $\left[ a,b\right] $ $\rightarrow \mathbb{R}^{+}$ be a multiplicative differentiable mapping on $\left[ a,b\right] $
with $a<b$. If $f^{\ast }$ is multiplicative integrable on $[a,b]$, then we
have the following identity for multiplicative integrals%
\begin{eqnarray*}
&&f\left( \tfrac{a+b}{2}\right) \left( \underset{a}{\overset{b}{\int }}%
f\left( u\right) ^{du}\right) ^{\frac{1}{a-b}} \\
&=&\left( \underset{0}{\overset{1}{\int }}\left( f^{\ast }\left( \left(
1-t\right) a+t\tfrac{a+b}{2}\right) ^{t}\right) ^{dt}\right) ^{\frac{b-a}{4}%
}\left( \underset{0}{\overset{1}{\int }}\left( f^{\ast }\left( \left(
1-t\right) \tfrac{a+b}{2}+tb\right) ^{t-1}\right) ^{dt}\right) ^{\frac{b-a}{4%
}}.
\end{eqnarray*}
\end{lemma}

\begin{proof}
Let%
\begin{equation*}
I_{1}=\left( \underset{0}{\overset{1}{\int }}\left( f^{\ast }\left( \left(
1-t\right) a+t\tfrac{a+b}{2}\right) ^{t}\right) ^{dt}\right) ^{\frac{b-a}{4}}
\end{equation*}%
and 
\begin{equation*}
I_{2}=\left( \underset{0}{\overset{1}{\int }}\left( f^{\ast }\left( \left(
1-t\right) \tfrac{a+b}{2}+tb\right) ^{t-1}\right) ^{dt}\right) ^{\frac{b-a}{4%
}}.
\end{equation*}%
Using the integration by parts for multiplicative integrals, from $I_{1}$ we
have%
\begin{eqnarray*}
I_{1} &=&\left( \underset{0}{\overset{1}{\int }}\left( f^{\ast }\left(
\left( 1-t\right) a+t\tfrac{a+b}{2}\right) ^{t}\right) ^{dt}\right) ^{\frac{%
b-a}{4}} \\
&=&\underset{0}{\overset{1}{\int }}\left( f^{\ast }\left( \left( 1-t\right)
a+t\tfrac{a+b}{2}\right) ^{\frac{b-a}{4}t}\right) ^{dt} \\
&=&\underset{0}{\overset{1}{\int }}\left( f^{\ast }\left( \left( 1-t\right)
a+t\tfrac{a+b}{2}\right) ^{\left( \frac{b-a}{2}\right) \frac{1}{2}t}\right)
^{dt} \\
&=&\tfrac{\left( f\left( \frac{a+b}{2}\right) \right) ^{\frac{1}{2}}}{1}.%
\tfrac{1}{\underset{0}{\overset{1}{\int }}\left( \left( f\left( \left(
1-t\right) \frac{a+b}{2}+tb\right) \right) ^{\frac{1}{2}}\right) ^{dt}} \\
&=&\left( f\left( \tfrac{a+b}{2}\right) \right) ^{\frac{1}{2}}\tfrac{1}{%
\left( \left( \underset{0}{\overset{1}{\int }}f\left( \left( 1-t\right) 
\frac{a+b}{2}+tb\right) ^{\frac{1}{2}}\right) ^{dt}\right) }=\left( f\left( 
\tfrac{a+b}{2}\right) \right) ^{\frac{1}{2}}\left( \underset{a}{\overset{%
\frac{a+b}{2}}{\int }}f\left( u\right) ^{du}\right) ^{\frac{1}{a-b}}.
\end{eqnarray*}%
Similarly, we have%
\begin{eqnarray*}
I_{2} &=&\left( \underset{0}{\overset{1}{\int }}\left( f^{\ast }\left(
\left( 1-t\right) \tfrac{a+b}{2}+tb\right) ^{t-1}\right) ^{dt}\right) ^{%
\frac{b-a}{4}} \\
&=&\underset{0}{\overset{1}{\int }}\left( f^{\ast }\left( \left( 1-t\right) 
\tfrac{a+b}{2}+tb\right) ^{\frac{b-a}{4}\left( t-1\right) }\right) ^{dt} \\
&=&\tfrac{1}{\left( f\left( \frac{a+b}{2}\right) \right) ^{-\frac{1}{2}}}.%
\tfrac{1}{\underset{0}{\overset{1}{\int }}\left( f\left( \left( 1-t\right) 
\frac{a+b}{2}+tb\right) \right) ^{dt}} \\
&=&\left( f\left( \tfrac{a+b}{2}\right) \right) ^{\frac{1}{2}}.\tfrac{1}{%
\left( \underset{0}{\overset{1}{\int }}\left( \left( f\left( \left(
1-t\right) \frac{a+b}{2}+tb\right) \right) ^{\frac{1}{2}}\right)
^{dt}\right) } \\
&=&\left( f\left( \tfrac{a+b}{2}\right) \right) ^{\frac{1}{2}}\left( 
\underset{0}{\overset{1}{\int }}f\left( \left( 1-t\right) \tfrac{a+b}{2}%
+tb\right) ^{dt}\right) ^{-\frac{1}{2}}=\left( f\left( \tfrac{a+b}{2}\right)
\right) ^{\frac{1}{2}}\left( \underset{\frac{a+b}{2}}{\overset{b}{\int }}%
f\left( u\right) ^{du}\right) ^{\frac{1}{a-b}}.
\end{eqnarray*}%
Multiplying above equalities we get the desired result. \\ The proof is
completed.
\end{proof}

\begin{lemma}\label{Lem3}
Let $f:$ $\left[ a,b\right] $ $\rightarrow \mathbb{R}^{+}$ be a multiplicative differentiable mapping on $\left[ a,b\right] $
with $a<b$. If $f^{\ast }$ is multiplicative integrable on $[a,b]$, then we
have the following identity for multiplicative integrals%
\begin{equation}
G\left( f\left( a\right) ,f\left( b\right) \right) \left( \underset{a}{%
\overset{b}{\int }}f\left( u\right) ^{du}\right) ^{\frac{1}{a-b}}=\left( 
\underset{0}{\overset{1}{\int }}\left( f^{\ast }\left( \left( 1-t\right)
a+tb\right) ^{\left( 2t-1\right) }\right) ^{dt}\right) ^{\frac{b-a}{2}},
\end{equation}%
where $G$ is the geometric mean i.e. $G\left( \alpha ,\beta \right) =\sqrt{%
\alpha \beta }$ for $\alpha ,\beta >0$.
\end{lemma}

\begin{proof}
Using integration by parts for multiplicative integrals, on the right side
of (3.1) we have%
\begin{eqnarray*}
I &=&\left( \underset{0}{\overset{1}{\int }}\left( f^{\ast }\left( \left(
1-t\right) a+tb\right) ^{\left( 2t-1\right) }\right) ^{dt}\right) ^{\frac{b-a%
}{2}} \\
&=&\underset{0}{\overset{1}{\int }}\left( f^{\ast }\left( \left( 1-t\right)
a+tb\right) ^{\frac{b-a}{2}\left( 2t-1\right) }\right) ^{dt} \\
&=&\underset{0}{\overset{1}{\int }}\left( f^{\ast }\left( \left( 1-t\right)
a+tb\right) ^{\left( b-a\right) \left( t-\frac{1}{2}\right) }\right) ^{dt} \\
&=&\tfrac{\left( f\left( b\right) \right) ^{\frac{1}{2}}}{\left( f\left(
a\right) \right) ^{-\frac{1}{2}}}.\tfrac{1}{\underset{0}{\overset{1}{\int }}%
\left( f\left( \left( 1-t\right) a+tb\right) \right) ^{dt}} \\
&=&\sqrt{f\left( a\right) f\left( b\right) }\tfrac{1}{\left( \underset{0}{%
\overset{1}{\int }}f\left( \left( 1-t\right) a+tb\right) ^{dt}\right) }%
=G\left( f\left( a\right) ,f\left( b\right) \right) \left( \underset{a}{%
\overset{b}{\int }}f\left( u\right) ^{du}\right) ^{\frac{1}{a-b}}.
\end{eqnarray*}%
The proof is completed.
\end{proof}

\begin{theorem} \label{Thm1}
Let $f:$ $\left[ a,b\right] $ $\rightarrow \mathbb{R}^{+}$ be a multiplicative differentiable mapping on $\left[ a,b\right] $
with $a<b$. If $f^{\ast }$ is multiplicative convex on $[a,b]$, then we have 
\begin{equation*}
\left\vert f\left( \tfrac{a+b}{2}\right) \left( \underset{a}{\overset{b}{%
\int }}f\left( u\right) ^{du}\right) ^{\frac{1}{a-b}}\right\vert \leq \left(
\left( f^{\ast }\left( a\right) \right) \left( f^{\ast }\left( \tfrac{a+b}{2}%
\right) \right) ^{4}\left( f^{\ast }\left( b\right) \right) \right) ^{\frac{%
b-a}{24}}.
\end{equation*}
\end{theorem}
\begin{proof}
From Lemma (\ref{lem2}), properties of multiplicative integral and the multiplicative convexity of $f^{\ast }$, we have
\begin{eqnarray*}
&&\left\vert f\left( \tfrac{a+b}{2}\right) \left( \underset{a}{\overset{b}{%
\int }}f\left( u\right) ^{du}\right) ^{\frac{1}{a-b}}\right\vert \\
&\leq &\left( \exp \tfrac{b-a}{4}\underset{0}{\overset{1}{\int }}\left\vert
\ln \left( f^{\ast }\left( \left( 1-t\right) a+t\tfrac{a+b}{2}\right)
^{t}\right) \right\vert dt\right) \\
&&\times \left( \exp \tfrac{b-a}{4}\underset{0}{\overset{1}{\int }}%
\left\vert \ln \left( f^{\ast }\left( \left( 1-t\right) \tfrac{a+b}{2}%
+tb\right) ^{t-1}\right) \right\vert dt\right) \\
&=&\left( \exp \tfrac{b-a}{4}\underset{0}{\overset{1}{\int }}t\left\vert
\ln \left( f^{\ast }\left( \left( 1-t\right) a+t\tfrac{a+b}{2}\right)
^{t}\right) \right\vert dt\right) \\
&&\times \left( \exp \tfrac{b-a}{4}\underset{0}{\overset{1}{\int }}\left(
1-t\right) \left\vert \ln \left( f^{\ast }\left( \left( 1-t\right) \tfrac{a+b%
}{2}+tb\right) \right) \right\vert dt\right) \\
&\leq &\left( \exp \tfrac{b-a}{4}\underset{0}{\overset{1}{\int }}%
t\left\vert \ln \left( f^{\ast }\left( a\right) \right) ^{\left( 1-t\right)
}\left( f^{\ast }\left( \tfrac{a+b}{2}\right) \right) ^{t}\right\vert
dt\right) \\
&&\times \left( \exp \tfrac{b-a}{4}\underset{0}{\overset{1}{\int }}\left(
1-t\right) \left\vert \ln \left( f^{\ast }\left( \tfrac{a+b}{2}\right)
\right) ^{\left( 1-t\right) }\left( f^{\ast }\left( b\right) \right)
^{t}\right\vert dt\right) \\
&=&\left( \exp \tfrac{b-a}{4}\underset{0}{\overset{1}{\int }}\left( t\left(
1-t\right) \ln \left( f^{\ast }\left( a\right) \right) +t^{2}\ln \left(
f^{\ast }\left( \tfrac{a+b}{2}\right) \right) \right) dt\right) \\
&&\times \left( \exp \tfrac{b-a}{4}\underset{0}{\overset{1}{\int }}\left(
\left( 1-t\right) ^{2}\ln \left( f^{\ast }\left( \tfrac{a+b}{2}\right)
\right) +t\left( 1-t\right) \ln \left( f^{\ast }\left( b\right) \right)
\right) dt\right) \\
&=&\left( \exp \tfrac{b-a}{4}\left( \ln f^{\ast }\left( a\right) \underset{0}%
{\overset{1}{\int }}t\left( 1-t\right) dt+\ln f^{\ast }\left( \tfrac{a+b}{2}%
\right) \underset{0}{\overset{1}{\int }}t^{2}dt\right) \right) \\
&&\times \left( \exp \tfrac{b-a}{4}\left( \ln f^{\ast }\left( \tfrac{a+b}{2}%
\right) \underset{0}{\overset{1}{\int }}\left( 1-t\right) ^{2}dt+\ln
f^{\ast }\left( b\right) \underset{0}{\overset{1}{\int }}t\left( 1-t\right)
dt\right) \right) \\
&=&\left( \exp \tfrac{b-a}{4}\left( \tfrac{1}{6}\ln f^{\ast }\left( a\right)
+\tfrac{1}{3}\ln f^{\ast }\left( \tfrac{a+b}{2}\right) \right) \right) \\
&&\times \left( \exp \tfrac{b-a}{4}\left( \tfrac{1}{3}\ln f^{\ast }\left( 
\tfrac{a+b}{2}\right) +\tfrac{1}{6}\ln f^{\ast }\left( b\right) \right)
\right) 
\end{eqnarray*}
Then, we obtain
\begin{eqnarray*}
\left\vert f\left( \tfrac{a+b}{2}\right) \left( \underset{a}{\overset{b}{%
\int }}f\left( u\right) ^{du}\right) ^{\frac{1}{a-b}}\right\vert 
& \leq & \left( \left( f^{\ast }\left( a\right) \right) ^{\frac{1}{6}}\left(
f^{\ast }\left( \tfrac{a+b}{2}\right) \right) ^{\frac{2}{3}}\left( f^{\ast
}\left( b\right) \right) ^{\frac{1}{6}}\right) ^{\frac{b-a}{4}} \\
&=&\left( \left( f^{\ast }\left( a\right) \right) \left( f^{\ast }\left( 
\tfrac{a+b}{2}\right) \right) ^{4}\left( f^{\ast }\left( b\right) \right)
\right) ^{\frac{b-a}{24}}.
\end{eqnarray*}
 The proof is completed.
\end{proof}
\begin{corollary}
In Theorem (\ref{Thm1}), If we assume that $f^{\ast }\leq M$, we get%
\begin{equation*}
\left\vert f\left( \tfrac{a+b}{2}\right) \left( \underset{a}{\overset{b}{%
\int }}f\left( u\right) ^{du}\right) ^{\frac{1}{a-b}}\right\vert \leq M^{%
\frac{b-a}{4}}.
\end{equation*}
\end{corollary}

\begin{corollary}\label{corl3}
In Theorem (\ref{Thm1}), using the multiplicative convexity of $f^{\ast }$i.e.\\ $$f^{\ast
}\left( \tfrac{a+b}{2}\right) \leq \sqrt{f^{\ast }\left( a\right) f^{\ast
}\left( b\right) },$$ we obtain%
\begin{equation*}
\left\vert f\left( \tfrac{a+b}{2}\right) \left( \underset{a}{\overset{b}{%
\int }}f\left( u\right) ^{du}\right) ^{\frac{1}{a-b}}\right\vert \leq \left(
f^{\ast }\left( a\right) f^{\ast }\left( b\right) \right) ^{\frac{b-a}{8}}.
\end{equation*}
\end{corollary}

\begin{theorem}\label{Thm6}
Let $f:$ $\left[ a,b\right] $ $\rightarrow \mathbb{R}^{+}$ be a multiplicative differentiable mapping on $\left[ a,b\right] $
with $a<b$. If $f^{\ast }$ is multiplicative convex on $[a,b]$, then we have 
\begin{equation*}
\left\vert G\left( f\left( a\right) ,f\left( b\right) \right) \left( 
\underset{a}{\overset{b}{\int }}f\left( u\right) ^{du}\right) ^{\frac{1}{a-b%
}}\right\vert \leq \left( f^{\ast }\left( a\right) f^{\ast }\left( b\right)
\right) ^{\tfrac{b-a}{8}},
\end{equation*}%
where $G$ is the geometric mean.
\end{theorem}
\begin{proof}
From Lemma (\ref{Lem3}), properties of multiplicative integral and the multiplicative convexity of $f^{\ast }$, we have
\begin{eqnarray*}
&&\left\vert \left( \left( f\left( a\right) \right) ^{2}\left( f\left( 
\tfrac{a+b}{2}\right) \right) ^{-1}f\left( b\right) ^{2}\right) ^{\frac{1}{3}%
}\left( \underset{a}{\overset{b}{\int }}f\left( u\right) ^{du}\right) ^{%
\frac{1}{a-b}}\right\vert \\
&\leq &\exp \tfrac{b-a}{2}\underset{0}{\overset{1}{\int }}\left( \ln
\left\vert f^{\ast }\left( \left( 1-t\right) a+tb\right) ^{\left(
2t-1\right) }\right\vert \right) dt \\
&=&\exp \tfrac{b-a}{2}\underset{0}{\overset{1}{\int }}\left( \left\vert
2t-1\right\vert \ln \left\vert f^{\ast }\left( \left( 1-t\right) a+tb\right)
\right\vert \right) dt \\
&\leq &\exp \tfrac{b-a}{2}\underset{0}{\overset{1}{\int }}\left( \left\vert
2t-1\right\vert \ln \left( f^{\ast }\left( a\right) \right) ^{\left(
1-t\right) }\left( f^{\ast }\left( b\right) \right) ^{t}\right) dt \\
&=&\exp \tfrac{b-a}{2}\underset{0}{\overset{1}{\int }}\left( \left\vert
2t-1\right\vert \left( \left( 1-t\right) \ln \left( f^{\ast }\left( a\right)
\right) +t\ln \left( f^{\ast }\left( b\right) \right) \right) \right) dt \\
&=&\exp \tfrac{b-a}{2}\left( \ln \left( f^{\ast }\left( a\right) \right) 
\underset{0}{\overset{1}{\int }}\left\vert 2t-1\right\vert \left(
1-t\right) dt+\ln \left( f^{\ast }\left( b\right) \right) \underset{0}{%
\overset{1}{\int }}\left\vert 2t-1\right\vert tdt\right) \\
&=&\left( \exp \tfrac{b-a}{2}\left( \ln \left( f^{\ast }\left( a\right)
\right) \underset{0}{\overset{1}{\int }}\left\vert 2t-1\right\vert \left(
1-t\right) dt\right) \right) \\
&&\times \left( \exp \tfrac{b-a}{2}\left( \ln \left( f^{\ast }\left(
b\right) \right) \underset{0}{\overset{1}{\int }}\left\vert 2t-1\right\vert
tdt\right) \right) \\
&=&\left( \exp \tfrac{b-a}{2}\left( \tfrac{1}{4}\ln \left( f^{\ast }\left(
a\right) \right) \right) \right) \left( \exp \tfrac{b-a}{2}\left( \tfrac{1}{4%
}\ln \left( f^{\ast }\left( b\right) \right) \right) \right) \\
&=&\left( f^{\ast }\left( a\right) f^{\ast }\left( b\right) \right) ^{\tfrac{%
b-a}{8}}.
\end{eqnarray*}%
The proof is completed.
\end{proof}

\begin{corollary}\label{corol1}
In Theorem (\ref{Thm6}), If we assume that $f^{\ast }\leq M$, we get%
\begin{equation*}
\left\vert G\left( f\left( a\right) ,f\left( b\right) \right) \left( 
\underset{a}{\overset{b}{\int }}f\left( u\right) ^{du}\right) ^{\frac{1}{a-b%
}}\right\vert \leq M^{\tfrac{b-a}{4}}.
\end{equation*}
\end{corollary}

\section{Applications to special means}

We shall consider the means for arbitrary real numbers $a,b$

The Arithmetic mean: $A\left( a,b\right) =\frac{a+b}{2}$.

The harmonic mean: $H\left( a,b\right) =\frac{2ab}{a+b}$, $a,b>0$.

The logarithmic means: $L\left( a,b\right) =\frac{b-a}{\ln b-\ln a}$, $%
a,b>0\ $and $a\neq b$.

The $p$-Logarithmic mean: $L_{p}\left( a,b\right) =\left( \frac{%
b^{p+1}-a^{p+1}}{\left( p+1\right) \left( b-a\right) }\right) ^{\frac{1}{p}}$%
, $a,b>0,a\neq b$ and $p\in 
\mathbb{R}
\backprime \left\{ -1,0\right\} $.

\begin{proposition}
Let $a,b\in 
\mathbb{R}
$ with $0<a<b$, then we have%
\begin{equation*}
e^{A^{p}\left( a,b\right) -L_{p}^{p}\left( a,b\right) }\leq \left(
e^{a^{p-1}+b^{p-1}}\right) ^{\frac{p\left( b-a\right) }{8}}.
\end{equation*}
\end{proposition}

\begin{proof}
The assertion follows from Corollary (\ref{corl3}), applied to the function $f\left(
t\right) =e^{t^{p}}$ with $p\geq 2$ whose $f^{\ast }\left( x\right)
=e^{pt^{p-1}}$ and $\left( \underset{a}{\overset{b}{\int }}f\left( u\right)
^{du}\right) ^{\frac{1}{a-b}}=\exp \left( -L_{p}^{p}\left( a,b\right)
\right) $.
\end{proof}

\begin{proposition}
Let $a,b\in \mathbb{R}$ with $0<a<b$, then we have%
\begin{equation*}
e^{H^{-1}\left( a,b\right) -L^{-1}\left( a,b\right) }\leq e^{-\tfrac{b-a}{%
4b^{2}}}.
\end{equation*}
\end{proposition}

\begin{proof}
The assertion follows from Corollary (\ref{corol1}) applied to the function $f\left(
t\right) =e^{\frac{1}{t}}$ whose $f^{\ast }\left( x\right) =e^{-\frac{1}{%
t^{2}}},$ $M=e^{-\frac{1}{b^{2}}}$ and $\left( \underset{a}{\overset{b}{\int 
}}f\left( u\right) ^{du}\right) ^{\frac{1}{a-b}}=\exp \left( -L^{-1}\left(
a,b\right) \right) $.
\end{proof}

\end{document}